\documentclass[12pt]{article}
\usepackage{amssymb,amsmath,amsfonts,amsthm,amscd,latexsym,verbatim,graphics,epsfig,indentfirst}
\usepackage{geometry} 
\def\pre{\mathrm{pre}}
\def\post{\mathrm{post}}
\def\RB{\mathrm{RB}}
\def\Var{\mathrm{Var}}
\def\Lie{\mathrm{Lie}}

\begin{document}

\begin{center}
{\Large Universal enveloping Lie Rota---Baxter algebras of
pre-Lie and post-Lie algebras}

V. Gubarev
\end{center}

\begin{abstract}
Universal enveloping Lie Rota---Baxter algebras
of pre-Lie and post-Lie algebras are constructed.
It is proved that the pairs of varieties
$(\RB\Lie,\mathrm{pre}\Lie)$ and $(\RB_\lambda\Lie$, $\mathrm{post}\Lie)$
are PBW-pairs and the variety of Lie Rota---Baxter algebras is not Schreier.

\end{abstract}

\section*{Introduction}

Pre-Lie algebras were firstly introduced by Cayley A. \cite{Cayley}
in the context of rooted tree algebras. In 1960s, pre-Lie algebras
appeared independently in affine geometry
(Vinberg E., \cite{Vinberg63}; Koszul J.-L. \cite{Koszul61}),
and ring theory (Gerstenhaber M., \cite{Gerst63}).
Arising from diverse areas, pre-Lie algebras are known under
different names like as Vinberg algebras, Koszul algebras,
left- or right-symmetric algebras (LSAs or RSAs), Gerstenhaber algebras.
Pre-Lie algebras satisfy an identity
$(x_1 x_2)x_3 - x_1(x_2 x_3) = (x_2 x_1) x_3 - x_2 (x_1 x_3)$.
See \cite{Bai,Burde06,Manchon} for surveys on pre-Lie algebras.

Post-Lie algebras were defined by Vallette B. in 2007 \cite{Vallette2007}
in the context of Loday algebras \cite{Dialg99,Trialg01}.
In last ten years, an amount of articles devoted to post-Lie algebras
in different areas is arisen \cite{Burde16,Fard,sl2}. A post-Lie algebra
is a vector space endowed with two bilinear products $[,]$ and $\cdot$
such that the bracket $[,]$ is Lie, and two identities are satisfied:
$$
(x \cdot y) \cdot z - x \cdot (y \cdot z)
- (y \cdot x) \cdot z + y \cdot (x \cdot z) = [y,x]\cdot z, \quad
x\cdot [y,z] = [x \cdot y,z] + [y,x\cdot z].
$$
If the bracket $[,]$ is zero, we have exactly a pre-Lie structure.

The varieties of pre- and post-Lie algebras play a key role
in the definition of any pre- and postalgebra
through black Manin operads product, see details in
\cite{BBGN2012,GubKol2014}.

In 2000, Golubchik I.Z., Sokolov V.V. \cite{GolubchikSokolov}
and Aguiar M. \cite{Aguiar00} independently noticed that starting with a 
Lie algebra $L$ over a field $\Bbbk$
endowed with a linear operator satisfying Rota---Baxter relation
\begin{equation}\label{eq:RB}
[R(x),R(y)] = R( [R(x),y] + [x,R(y)] + \lambda [x,y]), \quad x,y,\in A,\ \lambda\in \Bbbk,
\end{equation}
for $\lambda = 0$, we obtain a pre-Lie algebra structure
under the new operation $x\cdot y = [R(x),y]$.
In 2010 \cite{BaiGuoNi10}, Bai C., Guo L. and Ni X.
proved that a Lie algebra with a linear operator
satisfying \eqref{eq:RB} for $\lambda = 1$
is a post-Lie algebra with respect to the operations
$\cdot$ and $[,]$ where $x\cdot y = [R(x),y]$.

The operator satisfying \eqref{eq:RB} is called Rota---Baxter operator
of weight $\lambda$ (RB-operator, for short) and the algebra with RB-operator is called
Rota---Baxter algebra (RB-algebra). In 1960 \cite{Baxter60},
Baxter G. introduced these notions (in commutative case).
In 1980s, Belavin A.A., Drinfel'd V.G. \cite{BelaDrin82}
and Semenov-Tyan-Shanskii M.A. \cite{Semenov83}
stated a deep connection between Lie RB-algebras and Yang---Baxter equation.
About RB-algebras see a monograph of Guo L. \cite{Guo2011}.

In 2013 \cite{GubKol2013}, it was proved that any pre-Lie algebra
injectively embeds into its universal enveloping
Lie RB-algebra of zero weight and any post-Lie algebra
injectively embeds into its universal enveloping
Lie RB-algebra of nonzero weight.

Based on the last result, we have a problem:
To construct an universal enveloping Lie RB-algebra of pre- and post-Lie algebra.
Another connected problem is the following:
Whether the pairs of varieties $(\RB\Lie,\mathrm{pre}\Lie)$
and $(\RB_\lambda\Lie,\mathrm{post}\Lie)$, $\lambda\neq0$,
are PBW-pairs \cite{PBW}?
Here by $RB\Lie$ ($RB_\lambda\Lie$)
we mean the variety of Lie RB-algebras of (non)zero weight,
by $\mathrm{pre}\Lie$ and $\mathrm{post}\Lie$ --- the varieties of 
pre- and post-Lie algebras.

These two problems were stated in associative case were posted in \cite{Guo2011}
and were solved in \cite{Gub2017}.

The current work is devoted to the solution of the stated problems.
The solution is based on the construction of free Lie RB-algebra
obtained in \cite{Gub2016}.

In \S1, we give preliminaries on pre- and post-Lie algebras,
PBW-pairs of varieties, Lyndon---Shirshov words
and partially commutative Lyndon---Shirshov words.
Universal enveloping RB-algebras of pre-Lie (\S2)
and post-Lie algebras (\S3) are constructed.
In \S4, as corollaries we state that
the pairs of varieties $(\RB\Lie,\mathrm{pre}\Lie)$
and $(\RB_\lambda\Lie,\mathrm{post}\Lie)$ are PBW-pairs,
and the variety of Lie Rota---Baxter algebras is not a Schreier one.

\section{Preliminaries}

\subsection{Pre-Lie algebra}

Pre-Lie algebra is an algebra satisfying the pre-Lie identity:
\begin{equation}\label{pre-Lie}
(x_1 x_2)x_3 - x_1(x_2 x_3) = (x_2 x_1) x_3 - x_2 (x_1 x_3).
\end{equation}
By right-symmetric algebras (RSAs) one mean algebras
anti-isomorphic to pre-Lie algebras.

Any pre-Lie algebra is a Lie algebra under the commutator $[x,y] = xy - yx$.
Universal enveloping Lie algebra of a pre-Lie algebra was constructed by Segal D.
in 1994 \cite{Segal} (see also \cite{Oudom}).

{\bf Example 1} \cite{Burde06}.
Let $K = \Bbbk[x_1,\ldots,x_n]$ be a polynomial algebra with variables $x_1,\ldots,x_n$,
$\delta_i$ is a derivation of $K$ with respect to $x_i$.
Define on the space
$D_n = \bigg\{\sum\limits_{i=1}^n f_i\delta_i| f_i\in K\bigg\}$
of derivations of the algebra $K$ the following operation:
$f_i\delta_i \circ f_j\delta_j = f_i\delta_i (f_j)\delta_j$.
The space $D_n$ under the product $\circ$ is a pre-Lie algebra.

{\bf Example 2} \cite{ChapLiv}.
Consider a space $T$ of all rooted trees and define on $T$
the product $t_1*t_2$ as a sum of all trees obtained by
joining of the root of $t_1$ to every vertex of $t_2$.
The space $T$ under the product $*$ is free one-generated pre-Lie algebra.

In \cite{Burde06}, a good survey on pre-Lie algebras is given.
The problem of classifying of simple finite-dimensional pre-Lie algebras
is also considered.

An analog of the Magnus embedding theorem for right-symmetric algebras
is proved in \cite{KozUmir04}. The structure of universal
multiplicative enveloping algebras of free right-symmetric algebras
is studied in \cite{Kozybaev07}.

\subsection{Post-Lie algebra}

A post-Lie algebra is a vector space endowed with two bilinear products $[,]$ and $\cdot$,
the bracket $[,]$ is Lie, and the following identities are fulfilled:
\begin{gather}
(x \cdot y) \cdot z - x \cdot (y \cdot z)
- (y \cdot x) \cdot z + y \cdot (x \cdot z) = [y,x]\cdot z, \label{postLie1} \\
x\cdot [y,z] = [x \cdot y,z] + [y,x\cdot z]. \label{postLie2}
\end{gather}

In a context, universal enveloping Lie algebra
of a post-Lie algebra was found in \cite{Fard}.

\subsection{PBW-pair of varieties}

In 2014, Mikhalev A.A. and Shestakov I.P.
introduced the notion of a PBW-pair \cite{PBW};
it generalizes the relation between
the varieties of associative and Lie algebras
in the spirit of Poincar\'{e}---Birkhoff---Witt theorem.

Given varieties of algebras
$\mathcal{V}$ and $\mathcal{W}$, let
$\psi\colon \mathcal{V}\to \mathcal{W}$
be a such functor that maps an algebra $A\in \mathcal{V}$
to the algebra $\psi(A)\in \mathcal{W}$
preserving $A$ as vector space but changing the operations on $A$.
There exists left adjoint functor to the $\psi$
called universal enveloping algebra and denoted as $U(A)$.
Defining on $U(A)$ a natural ascending filtration,
we get associated graded algebra $\mathrm{gr}\,U(A)$.

Let $\mathrm{Ab}\,A$ denote the vector space $A$
with trivial multiplicative operations.

{\it Definition 1}.
A pair of varieties $(\mathcal{V},\mathcal{W})$
with the functor $\psi\colon \mathcal{V}\to \mathcal{W}$
is called PBW-pair if $\mathrm{gr}\,U(A)\cong U(\mathrm{Ab}\,A)$.

A variety $\mathcal{V}$ is called homogeneous if
the ideal of its identities is homogeneous, i.e.,
every homogeneous component of identities of $\mathcal{V}$
is itself an identity of $\mathcal{V}$.
A variety $\mathcal{V}$ is called Schreier
if every subalgebra of free $\mathcal{V}$-algebra is free.

{\bf Theorem 1} \cite{PBW}.
Let $(\mathcal{V},\mathcal{W})$ be a PBW-pair.
If $\mathcal{V}$ is a Schreier homogeneous variety, then
so is $\mathcal{W}$.

The varieties of pre- and post-Lie algebras,
and Lie RB-algebras of zero weight as well are homogeneous.

\subsection{Lyndon---Shirshov word}

Let $X$ be a well-ordered set with respect to an order
$<$, and let $X^*$ be the set of all associative words
in the alphabet $X$ (including the empty word which we denote by 1).
Extend the order to the set $X^*$ by induction on the word length as follows.
Put $u<1$ for every nonempty word $u$. Further,
$u < v$ for $u = x_i u'$, $v = x_jv'$, $x_i, x_j\in X$
if either $x_i < x_j$ or $x_i=x_j$, $u'< v'$.
In particular, the beginning of every word is greater than the whole word.

{\it Definition 2}.
A word $w\in X^*$ is called an associative Lyndon---Shirshov word if
for arbitrary nonempty $u$ and $v$ such that $w=uv$,
we have $w>vu$.

For example, a word $aabac$ is an associative Lyndon---Shirshov word when $a>b>c$.

Consider the set $X^+$ of all nonassociative words in $X$
(here we exclude the empty word from consideration), i.e.,
the words with all possible arrangements of parentheses.

{\it Definition 3}.
A nonassociative word $[u]\in X^+$ is called a nonassociative
Lyndon---Shirshov word (an $LS$-word, for short) provided that

(LS1) the associative word $u$ obtained from $[u]$
by eliminating all parentheses is an associative Lyndon---Shirshov word;

(LS2) if $[u] = [[u_1],[u_2]]$, then
$[u_1]$ and $[u_2]$ are $LS$-words, and $u_1> u_2$;

(LS3) if $[u_1] = [[u_{11}],[u_{12}]]$, then $u_2\geq u_{12}$.

These words appeared independently for the algebras and groups \cite{Shirshov1958,Lyndon1958}.
In \cite{Shirshov1958}, it was proved that the set of all
$LS$-words in the alphabet $X$ is a linear basis for a free Lie
algebra generated by $X$. Moreover, each associative Lyndon---Shirshov
word possesses the unique arrangement of parentheses which gives an LS-word.

\subsection{Partially commutative Lyndon---Shirshov word}

Partially commutative monoid was defined by Cartier P. and Foata D. in 1969 \cite{Traces69}.
Further, a lot of articles and monographs on partially commutative monoids, groups,
and algebras were wrtitten \cite{Kazachkov,Diekert,Duchamp,Duncan,Kim80}.

Let $G = \langle X, E\rangle$ be a (nonoriented) graph without
loops and multiple edges with the set of vertices $X$ and the set of edges $E$.

{\it Definition 4}.
A partially commutative Lie algebra $\Lie(G)$ with graph
$G = \langle X, E\rangle$ is a Lie algebra
with the set of generators $X$ and the defining relations
$[x_i,x_j] = 0$, where $(x_i, x_j)\in E$.

{\it Definition 5} \cite{Porosh}.
Define by induction partially commutative Lyndon---Shirshov words
($PCLS$-words) in an alphabet $X$ with a commutativity graph $G$:

(PCLS1) the elements of $X$ are $PCLS$-words;

(PCLS2) an $LS$-word $[u]$ of length greater than 1 in $X$ is a $PCLS$-word
provided that $[u]=([v],[w])$, where $[v]$ and $[w]$ are $PCLS$-words,
and the first letter of $[w]$ in $G$ is not connected
by an edge with at least one letter of $[v]$.

Poroshenko E. proved \cite{Porosh} that
the set of all $PCLS$-words defined by $X$ and a commutativity graph $G$
forms a linear base of a partially commutative Lie algebra $\Lie(G)$.
This result was obtained for a finite set $X$, but it may
be completely translated to a set of generators of arbitrary cardinality.

Given an alphabet $X$ with fixed commutativity graph $G$,
denote the sets of all $LS$- and $PCLS$-words by
$LS\langle X\rangle$ and $PCLS\langle X\rangle$ respectively.

\subsection{Free Lie RB-algebra}

Different linear bases (in different orders) of free Lie RB-algebra were constructed in
\cite{Gub2016,GubKol2016,Chen16}.

In \cite{Gub2016}, it was proved that due to the inductive construction
\begin{equation}\label{Z_infty}
Z_1 = LS\langle X\rangle,\quad
Z_2 = PCLS\langle X\cup R(Z_1)\rangle,\quad \ldots,\quad
Z_{n+1} = PCLS\langle X\cup R(Z_n)\rangle,\quad \ldots
\end{equation}
we obtain a set $Z_\infty = \cup Z_i$, a linear base of $\RB\Lie\langle X\rangle$.
Here we extend an order $<$ from $X$ to the set $\widetilde{X}$,
a union of $X$ and all elements of the form $R(w)$ ($R$-letter)
appeared in $Z_\infty$, as follows: $x<R(u)$ for all $x\in X$, $u\in \widetilde{X}$,
and $R(u)<R(v)$ if and only if $u<v$.
The commutativity graph with vertex set $\widetilde{X}$
is a clique on the set of all $R$-letters.

Let us call this base $Z_\infty$ of $\RB\Lie\langle X\rangle$ as a standard one.

Given a word $u$ from the standard base of $\RB\Lie\langle X\rangle$,
the number of appearances of the symbol $R$ in the notation of $u$
is called $R$-degree of the word $u$, denotation: $\deg_R(u)$.
Let us define a degree $\deg u$ of the word $u$ from the standard base
as the length of $u$ in the alphabet $\widetilde{X}$.
For example, given the word $u = [[R(R(x_1)),x_2],[R(x_3),x_4]]$,
we have $\deg_R(u) = 3$ and $\deg u = 4$.

{\bf Remark 1}.
Constructed standard base of $\RB\Lie\langle X\rangle$ does not depend
whether weight is zero or nonzero. The weight of RB-operator influences
only on the product of base elements.

\subsection{Embedding of pre- and post-Lie algebras in RB-algebras}

Further, unless otherwise specified,
RB-operator will mean RB-operator of zero weight.
We denote a free algebra of a variety $\Var$ generated by a set $X$ by
$\Var\langle X\rangle$, and a free RB-algebra of a variety $\Var$
and weight $\lambda$ respectively by
$\RB_\Var\Lie\langle X\rangle$. For short, denote
$\RB_0\Var\langle X\rangle$ by $\RB\Var\langle X\rangle$.

Given a Lie algebra $L$ with RB-operator $R$ of zero weight,
the space $L$ with respect to the operation
$x\cdot y = [R(x),y]$ is a pre-Lie algebra.

For any RB-operator $R$ of weight $\lambda\neq0$ on a Lie algebra $L$,
we have that an operator $R' = \frac{1}{\lambda}R$
is an RB-operator of unit weight.
The space $L$ with the operations $x\cdot y = [R(x),y]$ and $xy = [x,y]$
defined for $R'$ is a post-Lie algebra.
Denote the constructed pre- and post-Lie algebras as $L^{(R)}_{\lambda}$.
For short, we will denote $L^{(R)}_{0}$ as $L^{(R)}$.

Given a pre-Lie algebra $\langle L,\cdot\rangle$,
universal enveloping Lie RB-algebra $U$ of $L$
is an universal algebra in the class of all Lie RB-algebras
of zero weight such that there exists injective homomorphism from $L$ to $U^{(R)}$.
Analogously universal enveloping Lie RB-algebra of nonzero weight
of a post-Lie algebra is defined. The common denotation
of universal enveloping of pre- or post-Lie algebra: $U_{\RB}(L)$.

Let us write down a partilucar case of the result \cite{GubKol2013}:

{\bf Theorem 2} \cite{GubKol2013}.
Any pre-Lie (post-Lie) algebra could be embedded into its universal enveloping
Lie RB-algebra of (non)zero weight.

Based on Theorem 2, we have the natural question:
What does a linear base of universal enveloping RB-algebra
of a pre- or post-Lie algebra look like? Another problem is to determine
whether pairs $(\RB\Lie,\mathrm{pre}\Lie)$ and
$(\RB_\lambda\Lie$, $\mathrm{post}\Lie)$ are PBW-ones.
We will completely solve the problems.

In the article, the common method to construct universal enveloping is the following.
Let $X$ be a linear base of a pre-Lie algebra $L$.
We find a base of universal enveloping $U_{\RB}(L)$
as the special subset $E$ of the standard base of $\RB\Lie\langle X\rangle$
closed under the action of RB-operator.
By induction, we define a product $*$ on the linear span of $E$ and prove
that it satisfies Lie identities.
Finally, we state universality of the algebra $\Bbbk E$.

In the case of post-Lie algebra, as it was mentioned above, we will consider
universal enveloping Lie RB-algebra of unit weight.

\section{Universal enveloping Lie RB-algebra of a pre-Lie al\-gebra}

In the paragraph, we will construct universal enveloping Lie
RB-algebra $U_{\RB\Lie}(L)$ for arbitrary pre-Lie algebra $L$.

Let $X$ be a linear base of a pre-Lie algebra $\langle L,\cdot \rangle$.
Define a subset $E$ of a standard base of $\RB\Lie\langle X\rangle$ by induction:

1) $LS\langle X\rangle \subset E$;

2) given $u\in E$, define $R(u)\in E$;

3) given $u_1,\ldots,u_k\in E\setminus X$, $k\geq1$,
define that any (PCLS-)word of the degree not less than two
whose $R$-letters are exactly $R(u_1),\ldots,R(u_k)$ lies in $E$.

Due to the definition of the set $E$, we have three types of elements from $E$.

Let $D$ is a linear span of $E$, define the product $u*v$ on $D$
as it was done on the standard base of $\RB\Lie\langle X\rangle$ in \cite{Gub2016},
with the exception of two cases which are listed below.

Remember that product $u*v$ on $\RB\Lie\langle X\rangle$
was defined by induction on the parameters
$r = \deg_R(u)+\deg_R(v)$ and $q = \deg(u)+\deg(v)$.
We assumed that $u > v$, else we defined
$u*v = - v*u$ for $u\neq v$ and $u*u = 0$.
Proving bases of both inductions ($r = 0$, $q = 2$),
we fixed a set $F$ (in \cite{GubKol2016}, it was denoted as $E$)
of the letters from $\widetilde{X}$ containing at least in one of the words $u$ and $v$.
Thus, the set of words in alphabet $F$ of the degree $q$ is finite.
The third induction went on the decreasing of $\min(u,v)$ with respect to the order $<$.
Simultaneously with the definition of the product $*$,
it was proved by the third induction that
\begin{equation}\label{assump}
\begin{gathered}
u*v = \sum\limits_i \alpha_i w_i+f,\\
\deg f < q(u,v),\ \deg w_i = q(u,v),\ \alpha_i \in \Bbbk,
\ w_i = [w_{i1}, w_{i2}],\ w_{i2} \geq \min(u,v),
\end{gathered}
\end{equation}
all words $w_i$ are in alphabet $F$.
It was defined that $R(a)*R(b) := R(R(a)*b + a*R(b))$.

Now we write down the two cases in the definition of $u*v$ on $D$
which differ from the definition of $*$ on $\RB\Lie\langle X\rangle$:

Case 1: $R(x)*w$, $x\in X$, $w\in LS\langle X\rangle$.
Define such product by induction on $\deg w$.
For $\deg w = 1$, define $R(x)*y = x\cdot y$.
If $w = [w_1, w_2]$, then
$$
R(x)*[w_1,w_2] = (R(x)*w_1)*w_2 + w_1*(R(x)*w_2).
$$

Case 2: $R(x)*z$, $x\in X$, $z$ is a word from $E$ of type 3.
Define the product by induction on $\deg z$.
For $\deg z = 2$, define
$$
R(x)*[R(u),y] = [R(R(x)*u+x*R(u)),y] + [R(u),x\cdot y].
$$
If $z = [z_1, z_2]$, then
$$
R(x)*[z_1,z_2] = (R(x)*z_1)*z_2 + z_1*(R(x)*z_2).
$$

The product $u*v$ in the cases 1 and 2 is well-defined, as the products
$s_i = R(x)*w_i$ and $t_i = R(x)*z_i$, $i=1,2$, are yet defined by the inductions
on $\deg w$ and $\deg z$ respectively, and all products
$s_1*w_2$, $w_1*s_2$, $t_1*w_2$, $w_1*t_2$ are defined by the induction on
the sum of degrees.

Let $z$ be a word from $E$ of the type 1 or 3 equal to $z_1 z_2\ldots z_n$
by eliminating all operations $[,]$, where $z_i\in \widetilde{X}$.
It is easy to show that the product $*$ has the following property:
the product $R(x)*z$, $x\in X$, equals
$\sum\limits_{i=1}^n z(z_i\to R(x)*z_i)$,
where $z(z_i\to R(x)*z_i)$ denotes an expression
obtained from $z$ by replacing all sings $[,]$ on $*$
and replacing the $\widetilde{X}$-letter $z_i$ on the product $R(x)*z_i$.

{\bf Theorem 3}.
The space $D$ with respect to the operations $*$, $R$
is an universal enveloping Lie RB-algebra for a pre-Lie algebra $L$.

{\sc Proof}.
By the definition of $*$, the algebra $D$ is an anticommutative
enveloping RB-algebra for $L$. Let us prove that the product on $D$ satisfies the Jacobi identity:
$$
J(a,b,c) = (a*b)*c + (b*c)*a + (c*a)*b = 0.
$$
in the same way as it was proved in \cite{Gub2016}.
The proof in \cite{Gub2016} went by induction on
$r = \deg_R(a)+\deg_R(b)+\deg_R(c)$ and $q = \deg(a)+\deg(b)+\deg(c)$.
The proof of the base $r = 0$, of the base $q = 3$ for the triples
$a,b,c$ from $E$ of types 2--2--2 see in \cite{Gub2016}. The case $q = 3$
for the triples $a,b,c$ of types 2--1--1 follows from the definition $*$,
and for the triples $a,b,c$ of types 2--2--1 follows from \cite{Gub2016} and the equality
$$
J(R(x),R(y),z)
 = (x\cdot y)\cdot z - (y\cdot x)\cdot z - x\cdot (y\cdot z) + y\cdot (x\cdot z),
 \quad x,y,z\in X,
$$
satisfying in every pre-Lie algebra by \eqref{pre-Lie}.

After inductions bases on $r,q$, we fixed
a set $F$ (in \cite{GubKol2016}, it was denoted as $E$)
of the letters from $\widetilde{X}$ containing at least in one of the words $a,b,c$.
Thus, the set of words in alphabet $F$ of the degree $q$ is finite.

The inductive step on $q$ was proved by additional
induction among words in alphabet $F$ of the degree $q$
on decreasing of $c_0 = \min (a,b,c)$ and, if $c_0$ did not change,
on decreasing of $b_0 = \min (\{a,b,c\}\setminus\{c_0\})$.

By the definition of $*$ and the proof from \cite{Gub2016} it is enough to consider
only the cases when there are words $R(x)$, $x\in X$, among $a,b,c$.
We may assume that $a > b > c$.

{\sc Case 1}: $a = R(x)$, $x\in X$, $b = w_1$, $c = w_2$, $w_1,w_2\in LS\langle X\rangle$.
If $w_1*w_2 = [w_1,w_2]$, then $J(a,b,c)$ equals zero
when one computes $R(x)*[u_1,u_2]$ by the definition.
Let $w_1 = [w_{11},w_{12}]$ and $w_{12}>w_2$,
then specify the product $*$ on $LS\langle X\rangle$ as
как $w_1*w_2 = [[w_{11},w_2],w_{12}] + w_{11}*(w_{12}*w_2)$.
Such definition is well-defined by induction on decreasing
of the smallest from two multipliers of a $LS$-word, see \cite{Reut}.
Prove the Jacobi identity for the triple $a,b,c$
by the first induction on total degree and
the second decreasing induction on $\min(b,c)$.
For $\deg(b) + \deg (c) = 2$, it is obvious.
Notice that if $J(a,b_i,c_j) = 0$, $i=1,\ldots,k$, $j=1,\ldots,l$, then
$J(a,b,c) = 0$ for $b = \sum\limits_{i=1}^k \lambda_i b_i$, $c = \sum\limits_{j=1}^l \mu_j c_j$,
$\lambda_i,\mu_j\in\Bbbk$.
Hence, $J(R(x),b',c') = 0$ when yet $b',c'\not\in E$, but $b',c'$ are linearly expressed through
elements from $E$. Therefore, we can compute $R(x)*(w_1*w_2)$ by induction as
\begin{multline}\label{LieL:R(x)-w1-w2}\allowdisplaybreaks
R(x)*[[w_{11},w_2],w_{12}] + R(x)*(w_{11}*(w_{12}*w_2)) \\
= (R(x)*[w_{11},w_2])*w_{12} + [w_{11},w_2]*(R(x)*w_{12}) \\
+ (R(x)*w_{11})*(w_{12}*w_2) + w_{11}*(R(x)*(w_{12}*w_2)) \\
= ((R(x)*w_{11})*w_2)*w_{12} + (w_{11}*(R(x)*w_2))*w_{12} \\
+ [w_{11},w_2]*(R(x)*w_{12}) + (R(x)*w_{11})*(w_{12}*w_2) \\
+ w_{11}*((R(x)*w_{12})*w_2) + w_{11}*(w_{12}*(R(x)*w_2)).
\end{multline}
Subtracting from \eqref{LieL:R(x)-w1-w2} the following expression
\begin{multline*}
(R(x)*[w_{11},w_{12}])*w_2 + [w_{11},w_{12}]*(R(x)*w_2) \\
 = ((R(x)*w_{11})*w_{12})*w_2 + (w_{11}*(R(x)*w_{12}))*w_2 + [w_{11},w_{12}]*(R(x)*w_2)
\end{multline*}
and collecting summands by the groups with the product $R(x)*p$, $p\in\{w_{11},w_{12},w_2\}$,
we obtain zero. Thus, the Jacobi identity is satisfied.
Inductive applying of the Jacobi identity to the expressions $R(x)*p$ is correct,
as $R(x)*p$ equals a linear envelope of elements $E$ not containing $R$-letters.
The last one is zero by the Jacobi identity for $LS$-words.

{\sc Case 2}: $a = R(x)$, $b = R(y)$, $c = w\in LS\langle X\rangle$, $x,y\in X$.
Considering $w = [w_1,w_2]$ like an associative word, we have $w = r_1 r_2\ldots r_k$, $r_i\in X$.
Using the induction on $R$-degree, from one hand we compute
\begin{multline}\label{Lie1:R(x)-R(y)-w}
(R(x)*w)*R(y) = \sum\limits_{i=1}^k w(r_i\to x\cdot r_i)*R(y) \\
= - \sum\limits_{i=1}^k w(r_i\to y\cdot(x\cdot r_i))
- \sum\limits_{i,j=1,\,i\neq j}^k w(r_i\to x\cdot r_i; r_j\to y\cdot r_j);
\end{multline}
\vspace{-0.7cm}
\begin{multline}\label{Lie2:R(x)-R(y)-w}
(w*R(y))*R(x) = R(x)*\sum\limits_{i=1}^k w(r_i\to y\cdot r_i) \\
= \sum\limits_{i=1}^k w(r_i\to x\cdot(y\cdot r_i))
+ \sum\limits_{i,j=1,\,i\neq j}^k w(r_i\to y\cdot r_i; r_j\to x\cdot r_j).
\end{multline}
From another hand,
\begin{equation}\label{Lie3:R(x)-R(y)-w}
(R(y)*R(x))*w = R(y\cdot x - x\cdot y)*w
= \sum\limits_{i=1}^k w( r_i\to ((y\cdot x - x\cdot y)\cdot r_i) ).
\end{equation}
The sum of \eqref{Lie1:R(x)-R(y)-w}--\eqref{Lie3:R(x)-R(y)-w}
gives zero by \eqref{pre-Lie}.

{\sc Case 3}: $a = R(u)$, $b = R(x)$, $c = w\in LS\langle X\rangle$, $x\in X$, $u\not\in X$.
By the definition $R(u)*w = [R(u),w]$, and computing $[R(u),w]*R(x)$,
we prove the Jacobi identity.

{\sc Case 4}: $a = u_1$, $b = u_2$, $c = R(x)$, $x\in X$, $u_1,u_2$
are elements from $E$ of the type 3.
The proof is analogous to the proof of cases 1 and 2
in Theorem 1 \cite{Gub2016},
because there one do not need a specification of $c$.

{\sc Case 5}: $a = u$, $b = R(x)$, $c = w\in LS\langle X\rangle$,
$x\in X$, $u$ is an element from $E$ of the type 3.
If $u*w = [u,w]$, then calculating $R(x)*[u,w]$, we have done.
For $u_2>w$, where $u = [u_1,u_2]$,
we deal with it by analogous to the Case 1 induction on $c_0$.

{\sc Case 6}: $a = u$, $b = R(x)$, $c = R(y)$, $x,y\in X$,
$u$ is an element from $E$ of the type 3.
Computing all three summands from the Jacobi identity and applying the induction on the total degree,
we obtain a sum $\Sigma_1 + \Sigma_2 + \Sigma_3 + \Sigma_4$,
where $\Sigma_1$ is a result of the multiplication of a $R$-letter of $u$ on $R(x)$
and a letter from $X$ of $u$ on $R(y)$;
$\Sigma_2$ is a result of the multiplication of $R(x)$ and $R(y)$
on different $R$-letters of the word $u$;
$\Sigma_3$ is a result of the multiplication of $R(x)$ and $R(y)$
on the same $R$-letter of the word $u$.
Finally, $\Sigma_4$ is a result of the multiplication of $R(x)$ and $R(y)$
on letters from $X$ of the word $u$.
The sums $\Sigma_1$ and $\Sigma_2$ meet in the expressions
$(u*R(x))*R(y)$ and $(R(y)*u)*R(x)$ with the opposite signs and, therefore, give zero.
The sum $\Sigma_3$ is equal to zero by the same reason as the Jacobi identity is fulfilled 
in the case $q = 3$ and triple $a,b,c$ from $E$ of the types 2--2--2 (see \cite{Gub2016}).
The sum $\Sigma_4$ equals zero by \eqref{pre-Lie}.

{\sc Case 7}: $a = u$, $b = R(z)$, $c = R(x)$, $x\in X$, $z\not\in X$,
$u$ is an element from $E$ of the type 3.
The proof is analogous to the proof of cases 1 and 2 in Theorem 1 \cite{Gub2016},
because there one do not need a specification of $c$.

{\sc Case 8}: $a = R(z)$, $b = u$, $c = R(x)$, $x\in X$, $z\not\in X$,
$u$ is an element from $E$ of the type 3.
The proof is analogous to the proof of cases 5 and 6 in Theorem 1 \cite{Gub2016},
because there one do not need a specification of $c$.

Let us prove that the algebra $D$ is
exactly universal enveloping algebra for the pre-Lie algebra $L$,
i.e., is isomorphic to the algebra
$U_{\RB}(L) = \RB\Lie\langle X|x\cdot y = [R(x),y],\,x,y\in X\rangle$.
By the construction, the algebra $D$ is generated by $X$.
Therefore, $D$ is a homomorphic image of a homomorphism $\varphi$ from $U_{\RB}(L)$.

We will prove that all basic elements of $U_{\RB}(L)$ are linearly expressed by $D$,
this leads to nullity of kernel of $\varphi$ and $D\cong U_{\RB}(L)$.
Applying the identity
$$
[R(x),[u,v]] = [[R(x),u],v] + [u,[R(x),v]],\quad x\in X,
$$
we can prove by induction that the complement of the set $E$
in the base of $\RB\Lie(L)$ is linearly expressed via $E$.
This completes the proof.

\section{Universal enveloping Lie RB-algebra of a post-Lie algebra}

Let $X$ be a linear base of a post-Lie algebra $\langle L,\cdot,[,]\rangle$.
Define a subset $E$ of a standard base of $\RB\Lie\langle X\rangle$
by induction:

1) $X\subset E$;

2) given $u\in E$, define $R(u)\in E$;

3) given $u_1,\ldots,u_k\in E\setminus X$, $k\geq1$,
define that any (PCLS-)word of the degree not less than two
whose $R$-letters are exactly $R(u_1),\ldots,R(u_k)$ lies in $E$.

Let $D$ be a linear span of $E$, define the product $u*v$ on $D$
in analogous way as the product $*$ in the previous paragraph
except two cases.

At first, the product $x*y$ of $x,y\in X$ equals $x*y = [x,y]$.

At second, we define $R(u)*R(v) = R(R(u)*v + u*R(v) + u*v )$.

{\bf Theorem 4}.
The space $D$ with respect to the operations $*$, $R$
is an universal enveloping Lie RB-algebra of unit weight for $L$.

{\sc Proof}.
Analogous to the proof of Theorem 3.

\section{Corollaries}

{\bf Corollary 1}.
The pairs of varieties
$(\RB\Lie, \pre\Lie)$ and $(\RB_\lambda\Lie,\post\Lie)$ for $\lambda\neq0$
are PBW-pairs.

Notice that the varieties of Lie RB-algebras of zero weight
and ones of nonzero weight are simultaneously Schreier or not.
In \cite{Kozybaev07}, it was stated that the variety of pre-Lie algebras
is not a Schreirer one (see also \cite{Kozybaev08}).
Applying this result and Theorem 1, we obtain

{\bf Corollary 2}.
The variety of Lie RB-algebras is not a Schreirer variety.

{\bf Remark 2}.
Varieties of associative and commutative RB-algebras
are also not Schrei\-er. It is a particular case
of the common fact: given a not Schreier variety $\Var$,
we have that the variety of RB algebras of a variety $\Var$ is also not Schreier.
Indeed, let a subalgebra generated by a set $Y$ in $\Var\langle X\rangle$
(free algebra of a variety $\Var$ generated by $X$) be not free. 
Then a subalgebra generated by $Y$ in $\RB_\lambda\Var\langle X\rangle$
is also not free.

\section*{Acknowledgements}

The author is grateful to U. Umirbaev for the helpful references.

The research is supported by RSF (project N 14-21-00065).

\noindent
Gubarev Vsevolod \\
Sobolev Institute of Mathematics of the SB RAS \\
Acad. Koptyug ave., 4 \\
Novosibirsk State University  \\
Pirogova str., 2 \\
Novosibirsk, Russia, 630090\\
{\it E-mail: wsewolod89@gmail.com}

\end{document}